\documentclass[a4paper,12pt]{article}
\usepackage{amssymb,amsmath,amsfonts,pifont,xypic,amsthm,extarrows}

\begin{document}
\renewcommand{\theequation}{\arabic{section}.\arabic{equation}}
\numberwithin{equation}{section}

\newenvironment{zm}
               {\begin{sloppypar} \noindent{\bf Proof.}}
               {\hspace*{\fill} $\square$ \end{sloppypar}}

\newtheorem{atheorem}{\bf \temp}[section]
\newtheorem{dl}[atheorem]{Theorem}
\newtheorem{tl}[atheorem]{Corollary}
\newtheorem{yl}[atheorem]{Lemma}
\newtheorem{xz}[atheorem]{Property}
\newtheorem{mt}[atheorem]{Proposition}
\newtheorem{dy}[atheorem]{Definition}
\newtheorem{zs}[atheorem]{Remark}
\newtheorem{lz}[atheorem]{Example}
\newtheorem{cla}[atheorem]{Assertion}

\newtheorem*{zdl}{Main Theorem}

\title{Stability in graded rings associated with commutative augmented rings\footnote{Supported
by the NSFC (No. 11401155).}}
\author{\small Shan Chang \quad changshan@hfut.edu.cn}
\date{}
\maketitle
\vspace{-1cm}
\begin{center}
{\small School of Mathematics, Hefei University of Technology, Hefei 230009, China}
\end{center}

\begin{abstract}
Let $A$ be a commutative augmented ring and $I$ be its augmentation ideal.
This paper shows that the sequence $\{I^n/I^{n+1}\}$ becomes stationary up to isomorphism.
The result yields stability in the associated graded ring of $A$ along $I$. 
\end{abstract}

\section{Introduction}

A {\it commutative augmented ring} is defined to be a commutative ring $A$ together with
a homomorphism $\varepsilon: A\longrightarrow \mathbb{Z}$ which satisfies
\begin{itemize}
\item $A$ has an identity element and $\varepsilon$ preserves identity elements,
\item The underlying group of $A$ is a finitely generated free abelian group,
\item $I/I^2$ is torsion as an additive group, where $I=\ker \varepsilon$.
\end{itemize}
The homomorphism $\varepsilon$ is called the {\it augmentation map}. Its kernel $I$
is called the {\it augmentation ideal} of $A$.

The notion of commutative augmented ring arises from several natural rings associated
with finite groups. Let $G$ be a finite abelian group. The integral group ring $\mathbb{Z}G$
is a classical commutative augmented ring. Its augmentation map is induced by sending $g$
to $1$ for each element $g\in G$. Moreover, denote by $\mathcal{R}_H$ and $\Omega_H$
the representation ring and Burnside ring of a finite group $H$, respectively.
Then $\mathcal{R}_H$ and $\Omega_H$ are both commutative augmented rings.
Their augmentation maps are induced by the degree of representations and
the cardinality of fixed points of finite $H$-sets, respectively.
The details are presented in \cite{Pa}, \cite{CCT} and \cite{M}.

Let $A$ be a commutative augmented ring and $I$ be its augmentation ideal. Denote by $I^n$
and $Q_n(A)$ the $n$-th power of $I$ and the $n$-th consecutive quotient group $I^n/I^{n+1}$,
respectively. It is an interesting problem to determine the structures of $I^n$ and $Q_n(A)$
since they provide lots of information about $A$ itself. As a classical example,
it is well known that
\begin{equation}
\mbox{Tor}_1^A(A/I^n,A/I)\cong I^n/I^{n+1}=Q_n(A)
\end{equation}
for any positive integer $n$. The structure of augmentation ideals and their consective quotients for 
integeral group rings is one of most classical topic in group ring theory.  In particluar, Karpilovsky raised the
problem of determining the isomorphism type of the groups $Q_n(\mathbb{Z}G)$ in \cite{K}. This problem
has been well studied in \cite{BG}-\cite{CT}. The same problems for representation rings and Burnside rings
also been tackled in \cite{CCT}, \cite{C1}-\cite{CL} and \cite{WT1}-\cite{C4}, respectively .

There is a significant fact that all of known results show that the isomorphism
class of $Q_n$ does not depend on $n$ when $n$ is large enough. In particular, Bachman and
Grunenfelder showed the sequence $\{Q_n(\mathbb{Z}G)\}$ becomes stationary up to isomorphism
for any finite abelian group $G$ in \cite{BG}, the author, Chen and Tang obtained similar
results in \cite{CCT} for complex representation rings of dihedral groups and all finite abelian
groups, Wu and Tang proved in \cite{WT1} that, for any finite abelian group $G$, there exists a positive
integer $n_0$ such that $Q_n(\Omega_G)\cong Q_{n+1}(\Omega_G)$ for any $n\geqslant n_0$.
Motivated by this, we raise the following theorem. The proof is postponed till Section 4.
\begin{zdl}
Let $A$ be a commutative augmented ring. Then there exists a positive integer $n_0$
such that $Q_n(A)\cong Q_{n+1}(A)$ for any $n\geqslant n_0$.
\end{zdl}

A related problem of recent interest has been to investigate the augmentation ideals and
their consecutive quotients for integral group rings of finite non-abelian groups. The problem
has been tackled in \cite{ZT2}-\cite{ZN2}.

\section{Preliminaries}

In this section, we recall some classical results in commutative ring theory.
The following definition and lemmas were found in \cite{AM}. Through this section,
the word ring means a commutative ring with an identity element.
\begin{dy}
Let $\mathcal{C}$ be a class of $S$-module for some ring $S$,
and $\lambda$ be a function on $\mathcal{C}$ with value in $\mathbb{Z}$.
The function $\lambda$ is called additive if, for each short exact seuqence
\begin{equation}
0 \longrightarrow L' \longrightarrow L \longrightarrow L'' \longrightarrow 0
\end{equation}
in which all terms belong to $\mathcal{C}$, we have $\lambda(L')-\lambda(L)+\lambda(L'')=0$.
\end{dy}
\begin{yl}\label{th002}
Let $R=\bigoplus_{n=0}^{\infty} R_n$ be a graded ring. Then the following are equivalent.
\begin{itemize}
\item $R$ is a Noetherian ring.
\item $R_0$ is Noetherian and $R$ is finitely generated as an $R_0$-algebra.
\end{itemize}
\end{yl}
\begin{yl}\label{th003}
Let $M=\bigoplus_{n=0}^\infty M_n$ be a finitely generated graded $R$-module,
where $R=\bigoplus_{n=0}^\infty R_n$ is a Noetherian graded ring. Then
\begin{itemize}
\item $M_n$ is finitely generated as an $R_0$-module for each positive integer $n$,
\item If $R$ is generated by finite elements in $R_1$, then for all sufficiently large $n$,
$\lambda(M_n)$ is a polynomial in $n$ with rational coefficients, where
$\lambda$ is an additive function on the class of all finitely generated $R_0$-modules.
\end{itemize}
\end{yl}

\section{Necessary Tools}

This section proves some useful properties about commutative augmented rings and finite
abelian groups in this section. For convenience, we fix the following notations, where
$d$ and $N$ are positive integers, $p$ is a prime.
\begin{itemize}
\item Denote by $\mathbb{Z}_d$ the factor ring of $\mathbb{Z}$ modulo the ideal generated by $d$.
\item Denote by $\mathcal{A}$ the class of all finite abelian groups.
\item Denote by $\varepsilon_p(N)$ the exponent of $p$ in the unique factorization of $N$.
\item For any finite group $G$, denote the Sylow $p$-subgroup of $G$ by $G_{\{p\}}$
if $|G|$ is divided by $p$, otherwise, set $G_{\{p\}}$ to be the trivial subgroup.
\end{itemize}

The first lemma shows that we can construct Noetherian graded rings by
consecutive quotients of a commutative augmented ring naturally.
\begin{yl}\label{th004}
Let $A$ be a commutative augmented ring, and $I$ be its augmentation ideal. If $I/I^2$
is $d$-torsion as an additive group, then
\begin{itemize}
\item[\ding{172}] The graded ring $\mathcal{G}_A=\mathbb{Z}_d\oplus \Big(\bigoplus_{n=1}^\infty I^n/I^{n+1} \Big)$
is generated by finite elements in $I/I^2$ as an $\mathbb{Z}_d$-algebra,
\item[\ding{173}] $\mathcal{G}_A$ is Noetherian,
\item[\ding{174}] $\mathcal{P}_{A,m}=\bigoplus_{n=1}^\infty m(I^n/I^{n+1})$ is a finitely generated
graded $\mathcal{G}_A$-module for any positive integer $m$.
\end{itemize}
\end{yl}
\begin{zm}
\ding{172} Short calculation shows that $I^n/I^{n+1}$ is $d$-torsion as an abelian group
for each positive integer $n$. Thus $\mathcal{G}_A$ is a well-defined graded ring.
Note that $I$ is finitely generated as a free abelian group. Let $\{x_1,\ldots,x_r\}$
be a basis of $I$. Then $I/I^2$ is generated by $x_1+I^2,\ldots,x_r+I^2$ as a
$\mathbb{Z}_{d\,}$-module, which implies $\mathcal{G}_A=\mathbb{Z}_d[x_1+I^2,\ldots,x_r+I^2]$. \\
\ding{173} It is a direct corollary of \ding{172} because $\mathbb{Z}_d$ is obviously Noetherian. \\
\ding{174} $\mathcal{P}_{A,m}$ is finitely generated since it is actually an ideal of $\mathcal{G}_A$.
\end{zm}

The following two lemmas provide several additive functions on $\mathcal{A}$
and show the isomorphism class of a finite abelian group is determined by
the values of these additive functions on its subgroups.
\begin{yl}\label{th001}
For any fixed prime $p$, set
\begin{equation}
\lambda_p(G)=\varepsilon_p \big( |G| \big ), \quad G\in \mathcal{A}.
\end{equation}
Then $\lambda_p$ is an additive function on $\mathcal{A}$ with $\lambda_p(G)=\lambda_p(G_{\{p\}})$.
\end{yl}
\begin{zm}
For any short exact sequence
\begin{equation}
0 \longrightarrow G' \longrightarrow G \longrightarrow G'' \longrightarrow 0
\end{equation}
in $\mathcal{A}$, we have $|G|=|G'| \cdot |G''|$. From this it follows that
\begin{equation}
\varepsilon_p(|G|)=\varepsilon_p(|G'|)+\varepsilon_p(|G''|),
\end{equation}
hence $\lambda_p$ is additive. The last equality is implied by the fact that
$p$ does not divide $|G_{\{q\}}|$ for any prime $q\neq p$.
\end{zm}
\begin{yl}\label{th006}
Let $G$ and $G'$ be two finite abelian groups. If
\begin{equation}
\lambda_p(p^sG)=\lambda_p(p^sG')
\end{equation}
for any prime $p$ and any non-negative integer $s$, then $G\cong G'$.
\end{yl}
\begin{zm}
It is easy to verify that $(p^sG)_{\{p\}}=p^sG_{\{p\}}$ and $(p^sG')_{\{p\}}=p^sG'_{\{p\}}$.
So due to Lemma \ref{th001}, we can assume both $G$ and $G'$ are finite abelian $p$-groups.
Now suppose
\begin{equation}
G\cong \bigoplus_{k\geqslant 1} (\mathbb{Z}_{p^k})^{i_k}, \quad
G'\cong \bigoplus_{k\geqslant 1} (\mathbb{Z}_{p^k})^{j_k},
\end{equation}
where $i_k,j_k$ are non-negative integers and all but finitely many of them be zero.
Brief calculation shows that
\begin{equation}
p^sG\cong \bigoplus_{k\geqslant s+1} (\mathbb{Z}_{p^{k-s}})^{i_k}, \quad
p^sG'\cong \bigoplus_{k\geqslant s+1} (\mathbb{Z}_{p^{k-s}})^{j_k},
\end{equation}
hence
\begin{equation}
\lambda_p(p^sG)=\sum_{k\geqslant s+1} i_k(k-s), \quad \lambda_p(p^sG')=\sum_{k\geqslant s+1} j_s(k-s).
\end{equation}
From these it follows that
\begin{equation}\label{e001}
\sum_{k\geqslant s+1} (i_k-j_k)(k-s)=0.
\end{equation}
Denote by $\sigma_s$ the left side of (\ref{e001}). Then for any non-negative integer $s$, we have
\begin{equation}
\sigma_s-\sigma_{s+1}=\sum_{k\geqslant s+1}(i_k-j_k)=0.
\end{equation}
This implies $i_k=j_k$ for each positive integer $k$, the lemma is clear.
\end{zm}

\section{Main Theorem}

Now we prove the main theorem and apply it to integral group rings, representation rings
and Burnside rings.
\begin{zdl}
Let $A$ be a commutative augmented ring. Then there exists a positive integer $n_0$
such that $Q_n(A)\cong Q_{n+1}(A)$ for any $n\geqslant n_0$.
\end{zdl}
\begin{zm}
Fix an exponent $d$ of $I/I^2$. Then $\lambda_p$ is an additive function on the class
of all finite generated $\mathbb{Z}_{d\,}$-modules since these modules are finite abelian groups
with exponent $d$, where $p$ is a prime. Due to Lemma \ref{th003} and \ref{th004}, we get,
for each non-negative integer $s$, there is a polynomial $f_{p,s}$ in rational coefficients
such that
\[\lambda_p(p^sQ_n(A))=f_{p,s}(n)\]
when $n$ is large enough. To finish the proof, we need the following assertion.
\begin{cla}\label{th005}
For any positive integer $n$, we have $|Q_n(A)|\leqslant d^r$, where $r$ is the free rank of $I$.
\end{cla}
\begin{zm}
The proof of Lemma \ref{th004} shows $I^n/I^{n+1}$ is $d$-torsion for each positive integer $n$.
Then the assertion follows from the fact that $I^n$ is a free abelian group of the same rank as $I$.
\end{zm}
We return now to the proof of the theorem. By Assertion \ref{th005}, we get, for any positive
integer $n$,
\begin{equation}
\lambda_p(p^sQ_n(A))=\varepsilon_p \big( |p^sQ_n(A)| \big)
\leqslant \varepsilon_p \big( |Q_n(A)| \big)\leqslant \log_p(d^r),
\end{equation}
which implies the sequence $\{\lambda_p(p^sQ_n(A))\}_{n=1}^\infty$ is bounded. This forces the
polynomial $f_{p,s}$ be a constant. Hence $\lambda_p(p^sQ_n(A))$ is a constant for all
sufficiently large $n$. From this it follows that, for any fixed prime $p$ and non-negative integer
$s$, there is a positive integer $n_{p,s}$ such that
\begin{equation}
\lambda_p(p^sQ_n(A))=\lambda_p(p^sQ_{n+1}(A)), \quad n\geqslant n_{p,s}.
\end{equation}
Note that either $p\nmid d$ or $p^s>d^r$ implies $\lambda_p(p^sQ_n(A))=0$ for each positive integer $n$.
Thus we can set $n_{p,s}=1$ for all but finitely many $p$ and $s$. This implies $n_0=\max\{n_{p,s}\}$ is
well defined and
\begin{equation}
\lambda_p(p^sQ_n(A))=\lambda_p(p^sQ_{n+1}(A))
\end{equation}
holds for any prime $p$ and any non-negative integer $s$ whenever $n\geqslant n_0$.
Then the theorem follows from Lemma \ref{th006}.
\end{zm}

We have already mentioned that Bachman and Grunenfelder proved that, for any finite abelian
group $G$, the isomorphism class of $Q_n(\mathbb{Z}G)$ does not depend on $n$ when $n$ is large enough.
Here we regenerate this result independently as a direct corollary of Main Theorem.
\begin{tl}
For any finite abelian group $G$, there exists a positive integer $n_0$ such that
$Q_n(\mathbb{Z}G)\cong Q_{n+1}(\mathbb{Z}G)$ for any $n\geqslant n_0$.
\end{tl}

Furthermore, we get the following corollary by applying Main Theorem to representation
rings and Burnside rings.
\begin{tl}
Let $H$ be a finite group. Then there are two positive integers $n_1,n_2$ such that
\begin{itemize}
\item $Q_n(\mathcal{R}_H)\cong Q_{n+1}(\mathcal{R}_H)$ for any $n\geqslant n_1$,
\item $Q_n(\Omega_H)\cong Q_{n+1}(\Omega_H)$ for any $n\geqslant n_2$.
\end{itemize}
\end{tl}

\end{document}